\documentclass[12pt]{article}

\usepackage{graphicx}

\def\ra{\rightarrow}

\def\ss{\subseteq}

\def\e{\epsilon}

\def\Re{\hbox{\rm Re}\,}
\def\Im{\hbox{\rm Im}\,}

\def\dbar{\overline{\partial}}

 \def\HollowBox #1#2{{\dimen0=#1 \advance\dimen0 by -#2       
       \dimen1=#1 \advance\dimen1 by #2                       
        \vrule height #1 depth #2 width #2                    
        \vrule height 0pt depth #2 width #1                   
        \llap{\vrule height #1 depth -\dimen0 width \dimen1}% 
       \hskip -#2                                             
       \vrule height #1 depth #2 width #2}}                   
 \def\BoxOpTwo{\mathord{\HollowBox{6pt}{.4pt}}\;}             
\def\endpf{\hfill $\BoxOpTwo$}

\font\teneufm=eufm10
\font\seveneufm=eufm7
\font\fiveeufm=eufm5
\newfam\eufmfam
\textfont\eufmfam=\teneufm
\scriptfont\eufmfam=\seveneufm
\scriptscriptfont\eufmfam=\fiveeufm

\newfam\msbfam
\font\tenmsb=msbm10  scaled \magstep1 \textfont\msbfam=\tenmsb
\font\sevenmsb=msbm7 scaled \magstep1 \scriptfont\msbfam=\sevenmsb
\font\fivemsb=msbm5  scaled \magstep1 \scriptscriptfont\msbfam=\fivemsb
\def\Bbb{\fam\msbfam \tenmsb}

\def\RR{{\Bbb R}}
\def\CC{{\Bbb C}}

\def\NN{{\Bbb N}}

\newtheorem{theorem}{Theorem}

\newtheorem{proposition}[theorem]{Proposition}

\newtheorem{remark}[theorem]{Remark}

\newtheorem{definition}{Definition}
\newtheorem{example}[definition]{EXAMPLE}

\makeindex

\begin{document}

\begin{center}
\huge \bf A New Look at
\medskip \\
\huge \bf Convexity 
\medskip \\
\huge \bf and Pseudoconvexity
\end{center}

\begin{center}
\large Steven G. Krantz\footnote{Author supported in part
by a grant from the Dean of Graduate Studies at Washington University 
in St.\ Louis.  He thanks the American Institute of Mathematics for its
hospitality and support during this work.}
\end{center}
\vspace*{.15in}

\setcounter{section}{-1}

\section{Introduction}

Convexity is a classical idea.  Archimedes used a version of convexity in his
considerations of arc length.  Yet the idea was not formalized until 1934 in the
monograph of Bonneson and Fenchel [BOF].  
							   
The classical definition of convexity is this: An open domain $\Omega
\subseteq \RR^N$ is convex if, whenever $P, Q \in \Omega$, then the
segment $\overline{PQ}$ connecting $P$ to $Q$ lies in $\Omega$. We call
this the {\it synthetic} definition of convexity. It has the advantage of
being elementary and accessible (see [VAL]). The disadvantages are that it
is non-quantitative and non-analytic. It is of little use in situations of
mathematical analysis where it is most likely to arise.

The {\it analytic} definition of convexity is a bit more recondite.
Let $\Omega \ss \RR^N$ have $C^2$ boundary.  For us this means
that there exists a $C^2$ function $\rho$ defined in a neighborhood
$U$ of $\partial \Omega$ such that
$$
\Omega \cap U = \{x \in U: \rho(x) < 0\}
$$
and further that $\nabla \rho \ne 0$ on $\partial \Omega$.  We call $\rho$
a {\it defining function} for $\Omega$.  Let $P \in \partial \Omega$.
We say that a vector $w \in \RR^N$ is a {\it tangent vector} to $\partial \Omega$
at $P$, and we write $w \in T_P(\partial \Omega)$, if
$$
\sum_{j = 1}^N \frac{\partial \rho}{\partial x_j}(P) w_j = 0 \, .
$$
The domain $\Omega$ is said to be {\it analytically convex} at $P$ if
$$
\sum_{j, k = 1}^N \frac{\partial^2 \rho}{\partial x_j \partial x_k} (P) w_j w_k \geq 0  \eqno (*)
$$
for all $w \in T_P(\partial \Omega)$.

A moment's thought reveals that the condition $(*)$ simply mandates that the second
partial derivative of $\rho$ in the direction $w$ be nonnegative.  This is the
classical ``convex up'' condition from calculus.  This analytic definition of convexity
has the advantage that it can be localized to individual boundary points, and it
is {\it quantitative}.  It is a straightforward exercise (see [KRA1]) to see that
the analytic definition of convexity is equivalent to the synthetic definition
of convexity.

The notion of pseudoconvexity has a slightly different ontology.  Discovered by
E. E. Levi in the study of domains of holomorphy, this idea was first formulated
in its analytic form.  Let $\Omega \ss \CC^n$ have $C^2$ boundary.  Let
$\rho$ be a $C^2$ defining function for $\Omega$ as in our earlier discussion
of convexity.  Let $P \in \partial \Omega$.  We say that $\xi \in \CC^n$ is a 
{\it complex tangent vector} at $P$, and we write $\xi \in {\cal T}_P (\partial \Omega)$,
if
$$
\sum_{j = 1}^n \frac{\partial \rho}{\partial z_j}(P) \xi_j = 0 \, .
$$
The point $P$ is said to be a point of {\it Levi pseudoconvexity} if
$$
\sum_{j, k = 1}^n \frac{\partial^2 \rho}{\partial z_j \partial \overline{z}_k} (P) \xi_j \overline{\xi}_k \geq 0   \eqno (\star)
$$
for all $\xi \in {\cal T}_P(\partial \Omega)$.  

It is not a simple matter to give an elementary geometrical interpretation to the expression $(\star)$.
Part of the purpose of the present paper is to come to some basic geometric understanding of this notion of 
pseudoconvexity.

It is appropriate to record in passing a classical, alternative notion of pseudoconvexity.  Let $\Omega \ss \CC^n$ 
be any domain (smoothly bounded or not).  We say that $\Omega$ is {\it Hartogs pseudoconvex} if,
with $\delta_\Omega$ denoting the function of Euclidean distance to the boundary, we have
that $- \log \delta_\Omega$ is plurisubharmonic on $\Omega$.  It is known---see [KRA1]---that
a domain $\Omega$ with $C^2$ boundary is Levi pseudoconvex if and only if it is Hartogs pseudoconvex.

While the notion of Hartogs provides a sort of synthetic idea
of pseudoconvexity, it is not strictly analogous to the idea
that is used in classical convexity theory.  Convexity has played an ever
more prominent role in the function theory of several complex
variables in recent years (see [LEM] and [MCN1], [MCN2]).  Thus it
is worthwhile to be able to develop in further detail the
analogy between classical convexity theory and modern
pseudoconvexity theory.  That is our first purpose in the 
present paper.

A second purpose of this paper is to examine the concept of finite type.
Originating with the seminal paper [KOH] of Kohn, and later developed
by Catlin [CAT1], [CAT2] and D'Angelo [DAN], this idea has become
a central and influential artifact in complex function theory and
partial differential equations.  It is always worthwhile to find
new ways to understand finite type, and we explore some of these
in the present paper.

\section{Analytic Discs and Pseudoconvexity}

The results that we present here have a history.  Certainly they are related to
the classical {\it Kontinuit\"{a}tssatz}, for which see [KRA1].  But the proofs,
of necessity, are different.

Let $D \ss \CC$ be the unit disc. An {\it analytic disc} in
$\CC^n$ is a holomorphic mapping $\varphi: D \ra \CC^n$. A
{\it closed analytic disc} in $\CC^n$ is a continuous mapping
$\psi: \overline{D} \ra \CC^n$ such that $\psi \bigr |_D$ is
holomorphic. In practice we may refer to either of these
simply as an ``analytic disc''. The boundary of a closed
analytic disc is just $\psi(\partial D)$ whenever this
expression makes sense. It will frequently be convenient to
confuse the mapping $\varphi$ or $\psi$ with the image disc
$\varphi(D)$ or $\psi(\overline{D})$ (or the boundary map
$\psi(\partial D)$). We do so without further comment. The
center of an analytic disc is $\varphi(0)$ or $\psi(0)$.

In this paper we shall think of the boundary of a closed analytic disc as the complex-analytic
analogue of two points $P$ and $Q$ in the classical theory of convex sets. 
We shall think of the (image) analytic disc $\psi(\overline{D})$ as the 
complex-analytic analogue of the segment $\overline{PQ}$ that connects
$P$ and $Q$.  

Thus we should like to have a characterization of pseudoconvexity, in terms
of analytic discs, that is parallel to the synthetic characterization of convexity
in terms of segments.  It is the following.

\begin{proposition} \sl
Let $\Omega \ss \CC^n$ be a bounded domain with $C^2$ boundary.
Then $\Omega$ is Levi pseudoconvex in the classical sense if there
is a number $\delta_0 > 0$ so that, whenever $\psi: \overline{D} \rightarrow \CC^n$ 
is a closed analytic disc in $\CC^n$ with diameter less than $\delta_0$, 
and if $\psi(\partial D) \ss \partial \Omega$, then $\psi(\overline{D}) \ss \overline{\Omega}$.
\end{proposition}
{\bf Proof:}  Let \ dist \ denote Euclidean distance.  Choose $\epsilon_0 > 0$ so that 
$$
U_{\epsilon_0} \equiv \{z \in \CC^n: \hbox{dist}(z, \partial \Omega) < \epsilon_0\}
$$
is a tubular neighborhood of $\partial \Omega$ (see [HIR]).
Let $\delta_0 = \epsilon_0/100$.  Let $\psi$ be a closed analytic disc
as in the statement of the proposition.  It follows immediately from
the triangle inequality that the (image of the) closed analytic disc
lies entirely inside the tubular neighborhood $U_{\epsilon_0}$.\footnote{The purpose
of forcing the analytic disc to lie inside a tubular neighborhood (and to have small
diameter) is to guarantee
that the disc does not form the basis of a homology class in the boundary.}
Now there are two cases:
\begin{description}
\item[{\bf Some point of \boldmath $\psi(D)$ lies outside $\overline{\Omega}$}.]
In this case let $p_0 \equiv \psi(\zeta_0)$ be the point
of $\psi(D)$ that lies outside $\overline{\Omega}$ and furthest from $\partial \Omega$.  Let
$\nu$ be the unique normal vector from $\partial \Omega$ out to $p_0$.
Say that $\nu$ emanates from the base point $q_0 \in \partial \Omega$.
Then the domain
$$
\widehat{\Omega} \equiv \Omega - \nu = \{z - \nu: z \in \Omega\}
$$ 
has the property that the disc $\psi(D)$ is tangent to
$\partial \widehat{\Omega}$ at $q_0 - \nu$ and
the punctured disc $\psi(D) \setminus \{q_0 - \nu\}$ lies entirely
in $\widehat{\Omega}$.  But of course $\widehat{\Omega}$ is Levi pseudoconvex
with $C^2$ boundary.  So this last is impossible (see [KRA1]) by
the classical Kontinuit\"{a}tssatz.  We have eliminated this case.

\item[{\bf All points of \boldmath $\psi(D)$ lie in $\Omega$.}]
In this case $\psi(\overline{D}) \ss \overline{\Omega}$ and we are done.
\end{description}
\endpf
\smallskip \\

An argument similar to the one just presented, but even simpler, gives
the following result.  It is closer to the spirit of the classical
synthetic definition of convexity.

\begin{proposition} \sl
Let $\Omega \ss \CC^n$ be a bounded domain with $C^2$ boundary.
Then $\Omega$ is pseudoconvex in the classical sense if there
is a number $\delta_0 > 0$ so that, whenever $\psi: \overline{D} \rightarrow \CC^n$ 
is a closed analytic disc in $\CC^n$ with diameter less than $\delta_0$, 
and if $\psi(\partial D) \ss \Omega$, then $\psi(\overline{D}) \ss \Omega$.
\end{proposition}

Yet another variant is this:

\begin{proposition} \sl
Let $\Omega \ss \CC^n$ be a bounded domain with $C^2$ boundary.
Then $\Omega$ is pseudoconvex in the classical sense if there
is a number $\delta_0 > 0$ so that, whenever $\psi: \overline{D} \rightarrow \CC^n$ 
is a closed analytic disc in $\CC^n$ with diameter less than $\delta_0$, 
and if $\psi(\partial D) \ss \overline{\Omega}$, then $\psi(\overline{D}) \ss \overline{\Omega}$.
\end{proposition}

\section{The Concept of Finite Type}

The idea of finite type was first conceived in the paper [KOH] of Kohn.
Kohn's idea was to measure the complex-analytic flatness of a boundary
point of a domain in $\CC^2$; this was conceived as an obstruction
to subelliptic estimates for the $\overline{\partial}$-Neumann problem.

Later, Bloom and Graham [BLG] generalized Kohn's work to higher dimensions.
Perhaps more significantly, they isolated two very interesting definitions
of finite type and proved them to be equivalent.  We now briefly review
these two definitions.

\begin{definition} \rm
Let $\Omega \ss \CC^n$ be a domain with $C^\infty$ boundary.  Let $\rho$ be a smooth
defining function for $\Omega$.  Let $P \in \partial \Omega$.
Let $m$ be a positive integer.  We say that $P$ has {\it geometric type} at least $m$ if
there is a nonsingular analytic disc $\varphi: D \ra \CC^n$ such that $\varphi(0) = P$ and
$$
|\rho(\varphi(\zeta))| \leq C |\zeta|^m \, .  \eqno (\dagger)
$$
The greatest $m$ for which this is true is called the {\it type} of the point $P$.\footnote{And we
say that a disc satisfying condition $(\dagger)$ is {\it tangent to $\partial \Omega$ at $P$
to order $m$}.}  If there
is no greatest $m$ then the point $P$ is said to be of {\it infinite type}.
\end{definition}

Of course a point $P \in \partial \Omega$ has complex tangent space ${\cal T}_P(\partial \Omega)$ (see [KRA1]).
If $V$ is a small neighborhood of $P$ in $\partial \Omega$, then we may write down
a collection $L_1, \dots, L_{n-1}$ of tangent holomorphic vector fields on $V$ that are linearly independent
at each point of $V$.  A commutator (or Poisson bracket) $[L_j, L_k]$ or $[L_j, \overline{L}_k]$
or $[\overline{L}_j, \overline{L}_k]$ is called a {\it second-order commutator}.  If $M$ is
a $p^{\rm th}$-order commutator, then an expression of the form $[M, L_j]$ or $[M, \overline{L}_j]$ is called
a {\it $(p+1)^{\rm st}$-order commutator}.  For convenience, we refer to the individual
vector fields $L_1, \dots, L_{n-1}$ as {\it first-order commutators}.

\begin{definition} \rm
Let $\Omega \ss \CC^n$ be a domain with $C^\infty$ boundary.  Let $\rho$ be a smooth
defining function for $\Omega$.  Let $P \in \partial \Omega$.
Let $m$ be an integer exceeding 1.  We say that $P$ has {\it commutator type}
$m$ if any commutator $N$ of order $m-1$ or less satisfies
$$
\langle N, \partial \rho \rangle = 0
$$
but there is some commutator $N'$ of order $m$ that satisfies
$$
\langle N', \partial \rho \rangle \ne 0 \, .
$$
\end{definition}

The theorem of Bloom and Graham [BLG] says that, in $\CC^2$, a point $P \in \partial \Omega$
is of geometric type $m$ if and only if it is of commutator type $m$
(see [KRA1, pp.\ 463--464] for a quick and elegant proof).  In higher dimensions
this equivalence is still not fully understood, although there has been heartening recent
progress by Forn\ae ss and Lee [FOL].

John D'Angelo and David Catlin have demonstrated the importance of the concept
of finite type, both for function theory and for the study of the $\overline{\partial}$-Neumann
problem (see [DAN] and references therein).  See also the work of Baouendi, Ebenfelt, and Rothschild [BER]
for applications to the study of mappings.
It is worthwhile to be able to understand points of finite type from a variety of
different geometric points of view.

Our goal here is to understand the concept of finite type from the point of view
of analytic discs, analogous to our understanding of pseudoconvexity in the last section.
We continue to let \ dist \ denote Euclidean distance.  We also let
\ H-dist \ denote the Hausdorff distance on sets.  The result we are about to
present is certainly related to the work of Dwilewicz and Hill [DWH1], [DWH2].
These authors announce their results in all dimensions; but in the end
they only prove them in dimension two.  The results of the present
paper are valid in all dimensions.
%% are they??

\begin{proposition} \sl
Fix a domain $\Omega \ss \CC^n$ with smooth (that is, $C^\infty$) boundary.  
Let $P \in \partial \Omega$.  Fix an integer $m > 1$.  If $P$ has geometric
type $m$ then there is a sequence $\varphi_j: D \ra \Omega$ of analytic
discs satisfying
\begin{enumerate}
\item[{\bf (a)}]  $\varphi_j(0) \ra P$ as $j \ra \infty$.
\item[{\bf (b)}]  $\hbox{\rm diam}(\varphi_j(\overline{D})) \equiv \delta_j \ra 0$ as $j \ra \infty$.
\item[{\bf (c)}]  $\hbox{\rm H-dist}(\varphi_j(\overline{D}), \partial \Omega) \leq \delta_j^m$.
\end{enumerate}
\end{proposition}
{\bf Proof:}  Let $\varphi: D \ra \CC^n$ be an analytic disc that is
tangent to $\partial \Omega$ to order $m$ at $p$.  Let $\nu$ be the unit outward
normal vector to $\partial \Omega$ at $P$.  Then the discs
$$
\varphi_j = \varphi - \frac{1}{j} \nu
$$
(with an obvious slight adjustment to the domain of each disc)
will satisfy the three conclusions of the proposition.
\endpf
\smallskip \\

\def\btu{\bigtriangleup}

For the next result, which is the central one for this exposition, we need
to lay some groundwork.  First some background concepts and notation
(see [KRA3]).  If $g$ is a function on domain $U$ in $\RR^N$, $x \in U$, 
and $h \in \RR^N$ is small, then we let
$$
\btu_h g(x) = \btu^1_h g(x) \equiv g(x + h) - g(x)
$$
whenever this expression makes sense.  Iteratively, we set (for $j \in \NN$)
$$
\btu_h^j g(x) \equiv \btu_h (\btu_h^{j-1} g)(x) \, .
$$
So, for example, 
$$
\btu_h^2 g(x) = g(x + 2h) - 2 g(x + h) + g(x)
$$
and
$$
\btu_h^3 g(x) = g(x + 3h) - 3 g(x + 2h) + 3 g(x + h) - g(x) \, .
$$
We say ``whenever this expression makes sense'' because we must
ensure that $x, x+h$, etc., lie in $U$.

\begin{definition} \rm
Let $g$ be a continuous function on an open domain $U \subseteq \RR^N$.
Let $\alpha > 0$ and let $j$ be an integer that exceeds $\alpha$.  We
say that $g$ is {\it Lipschitz $\alpha$} on $U$, and write $g \in \Lambda_\alpha(U)$,
provided that $g$ is bounded on $U$ and there is a constant $C > 0$ such that
$$
\left | \btu_h^j g (x) \right | \leq C \cdot |h|^\alpha \, .
$$
\end{definition}

The paper [KRA3] discusses the equivalence of this definition with
other standard definitions of Lipschitz/H\"{o}lder functions.

\begin{proposition} \sl
Fix a domain $\Omega \ss \CC^2$ with smooth (that is, $C^\infty$) boundary.  
Let $P \in \partial \Omega$.  Fix an integer $m > 0$.  Assume that
there is a sequence $\varphi_j: D \ra \CC^n$ of analytic discs satisfying
\begin{enumerate}
\item[{\bf (a)}]  $\varphi_j(0) \ra P$ as $j \ra \infty$.
\item[{\bf (b)}]  $\hbox{\rm diam}(\varphi_j(\overline{D})) \equiv \delta_j \ra 0$ as $j \ra \infty$.
\item[{\bf (c)}]  $\hbox{\rm H-dist}(\varphi_j(\overline{D}), \partial \Omega) \leq \delta_j^m$.
\end{enumerate}
Then $P$ has geometric type at least  $m$.
\end{proposition}
{\bf Proof:}   It is a standard fact---see [KRA1, p.\ 463]---that there is a defining
function $\rho$ for $\Omega$ with the form
$$
\rho(z) = 2\Re z_2 + \mu(z_1) + {\cal O}(|z_1 z_2| + |z_2|^2) \, .
$$
In fact the first term of this expansion comes from normalizing the unit
outward normal at $(0,0) \in \partial \Omega$ to coincide with the real $z_2$
direction.  The second term comes from an application of E. Borel's theorem.
The remainder term comes from inspection of the Taylor expansion.

Now fix $\delta_j > 0$ as in the statement of the proposition
and let $\varphi_j: D \ra \Omega$ be an analytic disc satisfying
conditions {\bf (a)}, {\bf (b)}, {\bf (c)}. Certainly we may
suppose that this disc lies in a tubular neighborhood of the
boundary. We may as well suppose that this is a closed disc.
Therefore we may take a point $P = \varphi(\zeta_0)$ that is
furthest from the boundary. There is no loss of generality to
assume that $\zeta_0 = 0$ and we do so henceforth.

Let $\eta = \hbox{dist}(P, \partial \Omega)$ and let $\widetilde{P}$ be
the Euclidean projection of $P$ to the boundary of $\Omega$.  Let $\nu$
be the unit outward normal to $\partial \Omega$ at $\widetilde{P}$.
Now consider the analytic disc
$$
\widetilde{\varphi}(\zeta) = \varphi(\zeta) + \eta \nu \, .
$$
Assuming for the moment that $P$ is a nonsingular point of the analytic disc,\footnote{In
the case that the point is singular, we may perturb the analytic disc an arbitrarily
small amount in space to arrange that this extremal point {\it not} be singular.
After all, the set of singular points is discrete in the (image) analytic disc.}  we
now see that $\widetilde{\varphi}$ is tangent to $\partial\Omega$
at $\widetilde{P}$.  What is more, we may consider
$$
|\btu_h \rho \circ \widetilde{\varphi}(0)|  \eqno (*)
$$
for $h$ small (of size say $\delta_j/10$).  Condition {\bf (c)} in the statement of the proposition
tells us that the expression in $(*)$ is of size $\delta_j^m$.  The same
can be said for 
$$
|\btu_h^\ell \rho \circ \widetilde{\varphi}(0)|   \eqno (**)
$$
provided that $0 \leq \ell \leq m$.  Thus we see immediately (because
$\rho$ and $\mu$ are known in advance to be smooth functions)
that the function $\mu$ in the asymptotic expansion for $\rho$
vanishes to order $m$.  But this in turn says that the
analytic variety $\{z: z_2 = 0\}$ has order of contact $m$
with the boundary of $\Omega$ at 0.  So the point $0 \in \partial \Omega$ 
has type at least $m$.  That is what we wished to prove.
\endpf
\smallskip \\

\begin{remark} \rm Certainly our Proposition 5 implies the result of Dwilewicz
and Hill in [DWH1].  Our approach has the advantage that it does not
use delicate calculations involving the Bishop equation, and it generalizes
to higher dimensions.
\endpf
\end{remark}

The result of Proposition 5 is in complex dimension 2.  We wish
also to prove a result in complex dimension $n$ for any $n$.  We cannot
prove a sharp result at this time, but we can offer a useful characterization
of a finite type condition in terms of analytic discs.

\begin{proposition} \sl
Fix a domain $\Omega \ss \CC^n$ with smooth (that is, $C^\infty$) boundary.  
Let $P \in \partial \Omega$.  Fix an integer $m > 0$.  Assume that there
is a sequence of positive numbers $\delta_j \ra 0$ such that 
for any analytic disc $\varphi: D \ra \CC^n$ it holds that:
\begin{quote}
If $\hbox{dist}(\varphi(0), P) \approx \delta_j$ then
we have \ $\hbox{\rm H-dist}(\varphi(D), \partial \Omega) \leq \delta^m$.
\end{quote}
It follows there is a complex tangential direction $\tau$ along which $P$ has geometric type at least  $m$.
More precisely, there is a complex tangent vector $\tau$ and an analytic disc $\varphi: D \ra \CC^n$
such that
\begin{enumerate}
\item[{\bf (i)}]  $\varphi(0) = P$;
\item[{\bf (ii)}] $\varphi'(0) = c\tau$ for some complex constant $c$;
\item[{\bf (iii)}]  The disc $\varphi$ has order of contact at least $m$
with $\partial \Omega$ at $P$:
$$
|\rho \circ \varphi (P + \zeta\tau )| \leq C  |\zeta|^m \, .
$$
\end{enumerate}
\end{proposition}

\begin{remark} \rm It is natural to wonder under what circumstance the
hypothesis of Proposition 6 holds, near a given boundary point of
$\Omega$, for a sequence of $\delta$s tending to zero. It follows from
elementary contact geometry that the hypothesis will hold at {\it any}
smooth, pseudoconvex boundary point with $m = 2$. Thus any pseudoconvex
boundary point has type at least two, and that assertion is well known.
\endpf
\end{remark}

\begin{remark}  \rm It is worth noting that, in dimension $n = 2$, the hypothesis
of Proposition 6 is sufficient to imply the conclusion of Proposition 5.
In that case there is only one complex tangential direction, so the
conclusion of the proposition is considerably more elegant.
The proof is the same.
\endpf
\end{remark}

\noindent {\bf Proof of Proposition 6:}	  
Let $\delta_j \ra 0$ be chosen as in the statement of the proposition
and $\varphi_j$ be corresponding analytic discs.

We begin, just as in the proof
of Proposition 5, by writing a defining function for $\Omega$ as
$$
\rho(z) = 2 \Re z_n + \mu(z_1, \dots, z_{n-1}) + {\cal O}(|z_1 z_n| + |z_2 z_n| + \cdots |z_{n-1} z_n| + |z_n|^2) \, .  \eqno (\star)
$$
This, again, can be done by a normalization of coordinates (for the first term) and an application
of E. Borel's theorem (for the second term).  Arguing as before (and using
the notation of the proof of Proposition 5), we see that
$$
|\btu_h^\ell \rho \circ \widetilde{\varphi_j}(0)| \leq C |h|^m
$$
for $|h| \approx \delta$ and $0 \leq \ell \leq m$.  We may select a subsequence $\varphi_{j_k}$ such that
the tangent vectors $\varphi'_{j_k}(0)/\|\varphi'_{j_k}(0)\|$ converge to a limit
vector $\tau$.  It follows, just because $\rho$ and $\mu$ are smooth functions,
that
$$
\left ( \frac{\partial}{\partial \tau} \right )^j \mu(P) = 0 \quad \hbox{for} \ 0 \leq j \leq m \, .
$$
By our normalization $(\star)$ of the defining function, this gives the conclusion of
the proposition iwth the variety having contact $m$ being a complex line in the complex tangent
plane at $P$.
\endpf
\smallskip \\

\section{Other Geometric Conditions Involving Analytic Discs}

Let $\Omega$ be a smoothly bounded domain in $\CC^n$.  Let $P \in \partial \Omega$.
We now consider the following definition.

\begin{definition} \rm
Suppose that $\psi: \overline{D} \rightarrow \overline{\Omega}$ is a
closed analytic disc.  Assume that whenever $\psi(\partial D) \ss \partial \Omega$
then $\psi(D) \cap \partial \Omega = \emptyset$.  Then we say that every point
of $\partial \Omega$ is {\it complex analytically extreme}.
\end{definition}

This definition is analogous to the classical notion of ``extreme point'' from
the theory of convex sets (see, for example, [VAL]).  It is {\it not} the
case that a domain satisfying the condition of this last definition must have
the property that every boundary point is finite type.  The example
$$
\Omega = \{(z_1, z_2) \in \CC^2: |z_1|^2 + 2 e^{-1/|z_2|^2} < 1 \}
$$
illustrates this observation.

An analytically extreme boundary point is {\it analytically isolated} in the traditional
sense of that terminology.  This idea is important for the study of boundary
orbit accumulation points of automorphism group actions (for which see
[GK1], [GK2], and the survey [ISK]).

\section{Some Examples, and Comparison with Harmonic Discs}

\begin{example} \rm
Let 
$$
A = \{\zeta \in \CC: 1/2 < |\zeta| < 2\}
$$
and set
$$
\Omega = A \times D \, .
$$
Then $\Omega$ is a pseudoconvex domain in $\CC^2$.  It is a standard result---see [KRA1]---that
$\Omega$ may be exhausted by an increasing union of smoothly bounded, strongly pseudoconvex
domains $\Omega_j$.  Thus we may choose $j$ so large that
$$
\hbox{\rm dist}(\partial \Omega_j, \partial \Omega) < 10^{-10} \, .
$$
Now it is the case that the analytic disc
\begin{eqnarray*}
\psi: D & \ra & \Omega_j \\
     \zeta & \mapsto & (\zeta, 0)
\end{eqnarray*}
has the property that $\partial \psi$ lies in $\Omega_j$ while the
entire disc does not.  
\endpf
\end{example}

This example does not contradict our characterization of pseudoconvexity with
analytic discs (Proposition 1) because in that proposition we take
the discs sufficiently small that they do not generate any nontrivial homology classes
(see particularly the footnote on the next page).

\begin{example} \rm
Let 
$$
A = \{\zeta \in \CC: 1/2 < |\zeta| < 2\}
$$
and set
$$
\widetilde{\Omega} = (A \times D) \cup \biggl ( D(0,2) \times \{\zeta \in D: \hbox{Re}\, \zeta < -3/4\} \biggr ) \, .
$$
The disc
\begin{eqnarray*}
\psi: D & \ra & \widetilde{\Omega} \\
     \zeta & \mapsto & (\zeta, 0)
\end{eqnarray*}
still has the property that $\partial \psi(D) \ss \widetilde{\Omega}$.  Yet
the full disc does not lie in $\widetilde{\Omega}$.  

However note that this $\Omega$ is {\it not} pseudoconvex.  Indeed
the smallest pseudoconvex domain that contains $\Omega$ is $D(0,2) \times D$.
Thus, even though the disc $\psi$ has boundary curve that is homotopic to a point,
the example is insignificant because the domain is not pseudoconvex.
\endpf
\end{example}

It is natural to wonder whether analytic discs are the right device to use to measure pseudoconvexity.
Perhaps harmonic discs---which we use in effect to recognize plurisubharmonic functions---are more
appropriate.  They would certainly be more flexible.  Here by a {\it harmonic disc} we mean the following.
Let $\eta: \{e^{i\theta}: 0 \leq \theta \leq 2\pi\} \ra \CC^n$ be a continuous function with $\eta(e^{i0}) = \eta(e^{i2\pi})$.
Now solve the Dirichlet problem with boundary data $\eta$ to obtain a harmonic function
$$
u: D \ra \CC^n
$$
that is continuous up to the boundary and has boundary function $\eta$.  Is there a characterization of pseudoconvex domain using harmonic
discs that is analogous to Proposition 1?  The answer is ``no'', as the following example illustrates.

\def\e{\epsilon}

\begin{example} \rm
Let 
$$
A = \{\zeta \in \CC: |\Re \zeta| < 3, |\Im \zeta| < 3\} \setminus \{\zeta \in \CC: |\Re \zeta| \leq 1,
    |\Im \zeta| \leq 1\}
$$
and define
$$
\Omega = D(0,3) \times A \, .
$$
Now consider the curve $\gamma$, for $\e > 0$ small, that is traced out by
$$
\gamma(t) = \left \{ \begin{array}{lcr}
            (2 - \e, 1 + \e + i(1 + \e)t) & \hbox{\ if \ } & 0 \leq t \leq 1 \ \\ [.1in]
	    (2 - \e, 2 + 2\e + i(1 + \e) - t(1 + \e)) & \hbox{\ if \ } & 1 < t \leq 2 \ \\  [.1in]
	    (2 - \e, i(1 + \e) - 2(1 + \e) + t(1 + \e)) ) & \hbox{\ if \ } & 2 < t \leq 3 \ \\ [.1in]
	    (2 - \e, 1 + \e + 4i(1 + \e) - ti(1 + \e)) & \hbox{\ if \ } & 3 < t \leq 4 \, .
		     \end{array}
	    \right.
$$
This is certainly a closed curve {\it in} $\Omega$.  Yet the harmonic
disc with this curve as boundary (i.e., the solution of the Dirichlet problem with
boundary data given by $\gamma$) certainly has in it the point which is the average
over the curve.  A simple calculation shows that that point is $(2 - \e, 3/4(1 + \e) + 3i/4(1 + \e)) \not \in \Omega$.

Thus, even though $\Omega$ is obviously pseudoconvex, the harmonic disc condition fails.
\endpf
\end{example}

\newpage

\noindent {\Large \sc References}
\smallskip \\

\begin{enumerate}

\item[{\bf [BER]}]  M. S. Baouendi, P. Ebenfelt, and L. Rothschild, {\it Real
Submanifolds in Complex Space and their Mappings}, Princeton University Press,
Princeton, NJ, 1999.

\item[{\bf [BLG]}] T. Bloom and I. Graham, A geometric characterization of
points of type $m$ on real submanifolds of $\CC^n,$ {\em J. Diff. Geom.}
12(1977), 171-182.

\item[{\bf [BOF]}]  T. Bonnesen and W. Fenchel, {\it Theorie der Konvexen K\"{o}rper},
Springer-Verlag, Berlin, 1934.

\item[{\bf [CAT1]}] D. Catlin, Necessary conditions for
subellipticity of the $\overline{\partial}-$Neumann problem,
{\em Ann. Math.} 117(1983), 147-172.

\item[{\bf [CAT2]}] D. Catlin, Subelliptic estimates for the
$\dbar$Neumann problem, {\em Ann. Math.} 126(1987), 131-192.
		
\item[{\bf [DAN]}]  J. P. D'Angelo, {\em Several Complex Variables and
the Geometry of Real Hypersurfaces}, CRC Press, Boca Raton, 1993.
				  
\item[{\bf [DIF]}] K. Diederich and J. E. Forn\ae ss, Pseudoconvex domains:
Bounded strictly plurisubharmonic exhaustion functions, {\em Invent.
Math.} 39(1977), 129-141.

\item[{\bf [DWH1]}]  R. Dwilewicz and C. D. Hill, An analytic disc approach to the
notion of type of points, {\it Indiana Univ. Math. J.} 41(1992), 713--739.

\item[{\bf [DWH2]}]  R. Dwilewicz and C. D. Hill, A characterization of harmonic
functions and points of finite and infinite type, {\it Indag.\ Math.} (N.S.) 4(1993),
39--50.

\item[{\bf [FOL]}]  J. E. Forn\ae ss and L. Lee, personal communication.

\item[{\bf [GK1]}] R. E. Greene and S. G. Krantz, Invariants of Bergman
geometry and results concerning the automorphism groups of domains in
$\CC^n,$ Proceedings of the 1989 Conference in Cetraro (D. Struppa, ed.),
1991.

\item[{\bf [GK1]}] R. E. Greene and S. G. Krantz, Geometric Foundations for
Analysis on Complex Domains, Proceedings of the 1994 Conference in Cetraro
(D. Struppa, ed.), 1995.

\item[{\bf [HIR]}]  M. Hirsch, {\it Differential Topology}, Springer-Verlag,
New York, 1976.

\item[{\bf [ISK]}] A. V. Isaev and S. G. Krantz, Domains with non-compact
automorphism group: A Survey, {\it Advances in Math.} 146(1999), 1--38.

\item[{\bf [KOH]}] J. J. Kohn, Boundary behavior of $\dbar$ on weakly
pseudoconvex manifolds of dimension two, {\em J. Diff. Geom.} 6(1972),
523-542.

\item[{\bf [KRA1]}]  S. G. Krantz, {\it Function Theory of Several Complex Variables},
$2^{\rm nd}$ ed., American Mathematical Society, Providence, RI, 2001.

\item[{\bf [KRA2]}]  S. G. Krantz, Characterizations of various domains of
holomorphy via $\dbar-$ estimates and applications to a problem of Kohn,
{\it Illinois J. Math.} 23(1979), 267-285.

\item[{\bf [KRA3]}]  S. G. Krantz, Lipschitz spaces, smoothness of functions, and approximation
theory,  {\it Expositiones Math.} 3(1983), 193-260.

\item[{\bf [LEE]}]  L. Lee, thesis, Washington University, 2007.

\item[{\bf [LEM]}] L. Lempert, La metrique \mbox{K}obayashi et
las representation des domains sur la boule, {\em Bull. Soc.
Math. France} 109(1981), 427-474.

\item[{\bf [MCN1]}]  J. McNeal, Convex domains of finite type, {\it J. Functional
Analysis} 108(1992), 361--373.

\item[{\bf [MCN2]}]  J. McNeal,  Estimates on the Bergman kernels of
convex domains, {\it Adv.\ in Math.} 109(1994), 108--139.

\item[{\bf [VAL]}]  F. Valentine, {\it Convex Sets},
McGraw-Hill, New York, 1964.
			    
\end{enumerate}

\end{document}